\documentclass[preprint,10pt,twocolumn,5p,authoryear]{elsarticle}

\usepackage{amsmath,amssymb,amsfonts}
\usepackage{dsfont}
\usepackage{mathrsfs}

\usepackage{xcolor}

\usepackage{bbm}
 
\newtheorem{thm}{Theorem}[section]

\newtheorem{rem}{Remark}[section]

\usepackage{stackengine}
\stackMath

\newproof{pf}{\textbf{Proof}}

\newcommand{\vertiii}[1]{{\vert\kern-0.25ex\vert\kern-0.25ex\vert #1 
    \vert\kern-0.25ex\vert\kern-0.25ex\vert}}

\definecolor{internationalkleinblue}{rgb}{0.0, 0.18, 0.65}

\usepackage{newfloat,algcompatible}

\usepackage[ruled,vlined]{algorithm2e}

\usepackage{caption}

\journal{Elsevier}

\begin{document}

\begin{frontmatter}

\title{A generalized Routh-Hurwitz criterion for the stability analysis of polynomials with complex coefficients: application to the PI-control of vibrating structures} 

\author[Namur]{Anthony Hastir\corref{cor1}}\ead{anthony.hastir@unamur.be} 
\author[Namur,Tokyo]{Riccardo Muolo}
\cortext[cor1]{Corresponding author}

\address[Namur]{Department of Mathematics and naXys, Namur Institute for Complex Systems, University of Namur, Rue Grafé 2, 5000 Namur, Belgium}   
\address[Tokyo]{Department of Systems and Control Engineering, Tokyo Institute of Technology, 2 Chome-12-1 Ookayama, Tokyo 152-8550, Japan}


\begin{abstract}       
The classical Routh-Hurwitz criterion is one of the most popular methods to study the stability of polynomials with real coefficients, given its simplicity and ductility. However, when moving to polynomials with complex coefficients, a generalization exists but it is rather cumbersome and not as easy to apply. In this paper, we make such generalization clear and understandable for a wider public. {\color{black} To this purpose, we have broken down the procedure in an algorithmic form, so that the method is easily accessible and ready to be applied.} After having explained the method, we demonstrate its use to determine the external stability of a system consisting of the interconnection between a rotating shaft and a PI-regulator. The extended Routh-Hurwitz criterion gives then necessary and sufficient conditions on the gains of the PI-regulator to achieve stabilization of the system together with regulation of the output. This illustrative example makes our formulation of the extended Routh-Hurwitz criterion ready to be used in several other applications. 
\end{abstract}

\end{frontmatter}


\section{Introduction}

\noindent
\label{sec:intro}
Let us consider the following $n$-th order polynomial
\begin{equation}\label{eq:poly1}
q(s)=s^{n}+\sum_{j=1}^{n}\left(a_j+i b_j\right) s^{n-j},
\end{equation}
where $i$ denotes the imaginary unit throughout this note.
We want to study its stability, i.e., all its roots need to have negative real part. If all the coefficients $b_j=0$ $\forall~j\in \{1,...,n\}$, meaning that we would be dealing with real coefficients, we would rely on the well-known Routh-Hurwitz criterion \cite{Routh1877,Hurwitz1895}, which provides a simple algorithm to verify the stability conditions. However, to study the stability of the polynomial \eqref{eq:poly1} is not trivial and the method to obtain the stability conditions is not as straightforward as its real analogous. In the literature there are some available tools, but they are often developed for specific cases and their applicability is not immediately clear for a general public. For instance, one finds the so-called Kharitonov's theorem, first introduced in \cite{Kharitonov_Real} for polynomials with real coefficients and then extended in \cite{Kharitonov_Complex} in the complex case. This theorem consists in determining the region where the roots of a polynomial are located based on the same conclusion obtained for several \textit{upper-} and \textit{lower-}polynomials. That is, polynomials with coefficients encapsulating the coefficients of the nominal polynomial, sometimes called \textit{interval polynomials}. One needs $4$ polynomials in the real case, while $8$ polynomials need to be used when the coefficients are complex. Few years later, these results have been revisited in \cite{Minnichelli_Desoer} and \cite{CallierDesoer_book} by taking an engineering oriented point of view, notably. The major drawback of such an approach is that the coefficients of the original polynomial are not used directly, making this method not systematic. \\
We found that the most general method is the one developed in \cite{Frank_1946} and then recalled in \cite{Xie}, which is the one we will discuss pedagogically in the following. To the best of our knowledge, it is the most natural and direct extension of the classical Routh-Hurwitz criterion in the complex case. However, as a simple counter-example on a polynomial of degree $3$  may highlight, the main result that is presented in \cite{Xie} is wrong when the degree of the considered polynomial is odd. The same mistake has been further repeated in \cite{YuYu} on a particular example in network theory. This constitutes one additional reason for this note to describe the method in a constructive and a rigorous way. With an eye on applications, such a criterion turns out to be useful for determining for instance the stability of a dynamical system whose dynamics exhibit complex coefficients. Such cases arise often in rotordynamics to describe the behavior of rotating shafts, as highlighted in e.g. \cite{Loewy_1970} and \cite{Barnett_1983}. Complex coefficients appear also in the dynamics of electrical networks as is described in \cite{Varrichio_Network_Circuits}{\color{black}, and in the dynamics on directed (asymmetric) networks \cite{hwang2005synchronization,Asllani1} and hypergraphs \cite{gallo2022synchronization,DeLellis_pinning}}. Applications to the analysis of spontaneous self-excitation in induction generators are also found in \cite{Bodson} where the authors recall the generalized criterion for polynomials with complex coefficients. This criterion has been applied to fractional systems with applications to population dynamic models in \cite{Bourafa}. The criterion discussed in this note should then also be useful for analyzing the stability or developing control methods for such systems and beyond. 
{\color{black} In Appendix \ref{app:comparison}, we will compare the newly developed method with a previous one consisting in doubling the degree of the polynomial and applying the classical Routh-Hurwitz criterion \cite{carletti}, showing the advantage of the former.}\\

The paper is organized as follows: the extension of the classical Routh-Hurwitz criterion is highlighted in Section \ref{sec:method} as an algorithm, in which we make the distinction between \textit{$n$ odd} and \textit{$n$ even}. {\color{black} Let us remark that the algorithm form has the sole purpose to break down the procedure in simple terms so that it is easily understandable by the broadest possible audience.} The case of $n=4$ is developed in Section \ref{sec:example}. 
An example build from rotordynamics is then considered in Section \ref{sec:application}: a Proportional-Integral (PI) action is applied to the system and the stability properties of the closed-loop system are analyzed using the results described in the previous Sections. {\color{black} We show that the stability conditions for the complex polynomial are straightforwardly obtained with the method developed in this paper.} Some conclusions are addressed in Section \ref{sec:ccl}. {\color{black} Let us stress again that} a lot of attention has been paid to describing the method in a pedagogical way in the form of an algorithm, which then could be straightforwardly implemented {\color{black} and put to use in applications, of control theory and beyond, where (in)stability of complex polynomials is needed.}


\section{General description of the method}
\noindent
\label{sec:method}

Let us again consider the $n$-th order polynomial given in \eqref{eq:poly1}. As a matter of notation, let $\mathbb{C}^-_\xi$ (resp. $\mathbb{C}^+_\xi$) denote the open subset $\{s\in\mathbb{C}, \mathfrak{Re}(s)<\xi\}$ (resp. $\{s\in\mathbb{C}, \mathfrak{Re}(s)>\xi\}$), $\xi\in\mathbb{R}$. We use also the notation $\vert A\vert$ for the determinant of the matrix $A$. The general algorithm that determines whether the roots of $q(s)$ in \eqref{eq:poly1} are in $\mathbb{C}_0^-$ is presented in Algorithm \ref{Alg:RH}.

\begin{algorithm*}
\begin{algorithmic}
\STATE{\textbf{Data: }The complex coefficients of the polynomial given in \eqref{eq:poly1}.}
\STATE{\textbf{Output: }Necessary and sufficient conditions for the roots of \eqref{eq:poly1} to be in $\mathbb{C}^-_0$.}
\STATE{\textbf{Initialization: } In the case where $n$ is even, construct a $2\times n$ matrix as follows:\\
\begin{center}
\begin{tabular}{ |p{2.8cm}|p{2.8cm}|p{0.5cm}|p{3.3cm}|p{2.8cm}| }
\hline
$a_1^{(1)} = a_1$ & $b_2^{(1)} = b_2$ & $\dots$ & $a_{n-1}^{(1)} = a_{n-1}$ & $b_n^{(1)} = b_n$\\
\hline
$b_1^{(1)} = a_{1}^{(1)}b_1 - b_2^{(1)}$ & $a_2^{(1)} = a_1^{(1)}a_2 - a_{3}^{(1)}$ & $\dots$ & $b_{n-1}^{(1)} = a_1^{(1)}b_{n-1}-b_n$ & $a_n^{(1)} = a_1^{(1)}a_n$\\
\hline
\end{tabular}
\end{center}
When $n$ is odd, construct the following $2\times n$ matrix:
\begin{center}
\begin{tabular}{ |p{2.8cm}|p{2.8cm}|p{0.5cm}|p{3.3cm}|p{2.8cm}| }
\hline
$a_1^{(1)} = a_1$ & $b_2^{(1)} = b_2$ & $\dots$ & $b_{n-1}^{(1)} = b_{n-1}$ & $a_n^{(1)} = a_n$\\
\hline
$b_1^{(1)} = a_{1}^{(1)}b_1 - b_2^{(1)}$ & $a_2^{(1)} = a_1^{(1)}a_2 - a_{3}^{(1)}$ & $\dots$ & $a_{n-1}^{(1)} = a_1^{(1)}a_{n-1}-a_n$ & $b_n^{(1)} = a_1^{(1)}b_n$\\
\hline
\end{tabular}
\end{center}
}
\FOR{$p=2,\dots,n-1$}
    \IF{($p$ is even and $n$ is even) or ($p$ is odd and $n$ is odd)}
        \STATE{Construct the following $2\times (n-(p-1))$ matrix:\\
        \begin{center}
        \begin{tabular}{ |p{0.8cm}|p{0.8cm}|p{0.5cm}|p{0.8cm}|p{0.8cm}| }
        \hline
        $a_p^{(p)}$ & $b_{p+1}^{(p)}$ & $\dots$ & $b_{n-1}^{(p)}$ & $a_n^{(p)}$\\
        \hline
        $b_p^{(p)}$ & $a_{p+1}^{(p)}$ & $\dots$ & $a_{n-1}^{(p)}$ & $b_n^{(p)}$\\
        \hline
        \end{tabular}
        \end{center}
        }
        where 
        \begin{itemize}
            \item The elements $a_k^{(p)}, k = p, p+2, p+4, \dots, n$ and $b_l^{(p)}, l = p+1, p+3, \dots, n-1$ of the first row are: $a_k^{(p)} = \left\vert\begin{matrix}a_{p-1}^{(p-1)} & -b_k^{(p-1)}\\ b_{p-1}^{(p-1)} & a_k^{(p-1)}\end{matrix}\right\vert$
            and $b_l^{(p)} = \left\vert\begin{matrix}a_{p-1}^{(p-1)} & a_l^{(p-1)}\\ b_{p-1}^{(p-1)} & b_l^{(p-1)}\end{matrix}\right\vert$
            \item The elements $a_k^{(p)}, k=p+1, p+3, \dots, n-1, b_l^{(p)}, l=p, p+2, p+4, \dots, n-2$ of the second row are: $a_k^{(p)} = -\left\vert\begin{matrix} a_{p-1}^{(p-1)}& a_k^{(p-1)}\\  a_p^{(p)}& a_{k+1}^{(p)}\end{matrix}\right\vert, b_l^{(p)} = -\left\vert\begin{matrix} a_{p-1}^{(p-1)}& b_l^{(p-1)}\\a_p^{(p)}  & b_{l+1}^{(p)}\end{matrix}\right\vert$ and $b_n^{(p)} = -\left\vert\begin{matrix} a_{p-1}^{(p-1)}& b_n^{(p-1)}\\a_p^{(p)}  & 0\end{matrix}\right\vert$  
        \end{itemize}
        \ENDIF
\IF{($p$ is odd and $n$ is even) or ($p$ is even and $n$ is odd)}
\STATE{Construct the following $2\times (n-(p-1))$ matrix:
\begin{center}
        \begin{tabular}{ |p{0.8cm}|p{0.8cm}|p{0.5cm}|p{0.8cm}|p{0.8cm}| }
        \hline
        $a_p^{(p)}$ & $b_{p+1}^{(p)}$ & $\dots$ & $a_{n-1}^{(p)}$ & $b_n^{(p)}$\\
        \hline
        $b_p^{(p)}$ & $a_{p+1}^{(p)}$ & $\dots$ & $b_{n-1}^{(p)}$ & $a_n^{(p)}$\\
        \hline
        \end{tabular}
        \end{center}
        }
        where 
        \begin{itemize}
            \item The elements $a_k^{(p)}, k = p, p+2, p+4, \dots, n-1$ and $b_l^{(p)}, l = p+1, p+3, \dots, n$ of the first row are: $a_k^{(p)} = \left\vert\begin{matrix}a_{p-1}^{(p-1)} & -b_k^{(p-1)}\\ b_{p-1}^{(p-1)} & a_k^{(p-1)}\end{matrix}\right\vert$
            and $b_l^{(p)} = \left\vert\begin{matrix}a_{p-1}^{(p-1)} & a_l^{(p-1)}\\ b_{p-1}^{(p-1)} & b_l^{(p-1)}\end{matrix}\right\vert$
            \item The elements $a_k^{(p)}, k=p+1, p+3, \dots, n-2$ and $b_l^{(p)}, l=p, p+2, p+4, \dots, n-1$ of the second row are: $a_k^{(p)} = -\left\vert\begin{matrix} a_{p-1}^{(p-1)}& a_k^{(p-1)}\\  a_p^{(p)}& a_{k+1}^{(p)}\end{matrix}\right\vert, b_l^{(p)} = -\left\vert\begin{matrix} a_{p-1}^{(p-1)}& b_l^{(p-1)}\\a_p^{(p)}  & b_{l+1}^{(p)}\end{matrix}\right\vert$ and $a_n^{(p)} = -\left\vert\begin{matrix} a_{p-1}^{(p-1)}& a_n^{(p-1)}\\a_p^{(p)}  & 0\end{matrix}\right\vert$  
        \end{itemize}
\ENDIF
\ENDFOR
\STATE{Compute the following coefficient
\begin{align*}
    a_n^{(n)} = \left\vert\begin{matrix} a_{n-1}^{(n-1)}& -b_n^{(n-1)}\\ b_{n-1}^{(n-1)}& a_n^{(n-1)} \end{matrix}\right\vert
\end{align*}
}
\STATE{\textbf{Necessary and sufficient conditions: } the roots of the polynomial \eqref{eq:poly1} are in $\mathbb{C}_0^-$ if and only if $a_k^{(k)} > 0$ for all $k=1,\dots,n$.}
\end{algorithmic}
\caption{Generalized Routh-Hurwitz criterion\label{Alg:RH}}
\end{algorithm*}

{\color{black}
\begin{rem}
Compared to Algorithm 1, the necessary and sufficient conditions obtained in \cite{Xie} for the polynomial (1) to be stable are stated as follows: $a_j^{(j)}>0, j = 1,\dots,\overline{n}$ with $\overline{n} = n$ is $n$ is even and $\overline{n} = n-1$ if $n$ is odd. Considering the particular case in which the imaginary parts of the coefficients are all $0$ and comparing this with the classical Routh-Hurwitz test for polynomials with real coefficients does not give the same conditions. By making a test on a general polynomial of degree $3$, we found out that $\overline{n}$ has to be equal to $n$ in the case of odd polynomials as well.
\end{rem}
}

{\color{black}
The proof that the proposed algorithm converges is given in \cite{Frank_1946}. The method on which it is based comes from \cite{Wall} and is centered around the representation of a certain ratio between two polynomials in terms of continued fractions with some properties. H.S. Wall in \cite{Wall} has been the first to prove the Routh criterion introduced in \cite{Hurwitz1895} for polynomials with real coefficients with a method based on continued fraction expansions. The complete proof of that result for polynomials with complex coefficients is not given in detail here but its main ingredients are recalled. Interested readers may have a look at \cite[Theorems 3.1 and 3.2]{Frank_1946}.\\
Let us first introduce the auxiliary polynomial $\mathfrak{p}(s)$ in the following way
\begin{equation}
    \mathfrak{p}(s) := \sum_{\overset{j=1}{j\text{ odd} }}^{n-1} a_j s^{n-j} + i\sum_{\overset{j=2}{j\text{ even}}}^{n-2} b_j s^{n-j}.
    \label{Aux_Poly}
\end{equation}
The main theorem is stated as follows, see \cite[Theorem 3.1]{Frank_1946}.
\begin{thm}\label{Thm_Cont_Frac}
    Let $p(s)$ and $\mathfrak{p}(s)$ be the polynomials given in (1) and \eqref{Aux_Poly}, respectively. The polynomial $p(s)$ is stable if and only if the ratio $\frac{\mathfrak{p}(s)}{p(s)}$ may be written as the following continued fraction
    \begin{equation}
        \displaystyle\frac{\mathfrak{p}(s)}{p(s)} = \frac{c_0}{s + c_0 + d_1 + \displaystyle\frac{c_1}{s + d_2 + \displaystyle\frac{c_2}{s + d_3 + \stackunder{}{\ddots\stackunder{}{\displaystyle{}+ \frac{c_{n-1}}{\displaystyle
        s+d_n}}}}}}
        \label{Continued_Frac}
    \end{equation}
    where $c_i, i=0, \dots, n-1$ are real and positive and $d_i, i=1, \dots, n$ are pure imaginary or zero.
\end{thm}
Then, \cite{Frank_1946} gives a characterization for the coefficients $c_i, i=0,\dots,n-1$ to be positive in terms of determinants that will be comparable to the numbers $a^{(k)}_k, k=1, \dots, n$, see \cite[Theorem 3.2]{Frank_1946}.
\begin{thm}\label{Thm_Determinants}
    The polynomial $p(s)$ is stable if and only if the following determinants are all positive.
    \begin{align*}
        &\Delta_1 = a_1,\\
        &\Delta_k =\\
        &\left\vert\begin{matrix}a_1 & a_3 & a_5 & \dots & a_{2k-1} & -b_2 & -b_4 & \dots & -b_{2k-2}\\
        1 & a_2 & a_4 & \dots & a_{2k-2} & -b_1 & -b_3 & \dots & -b_{2k-3}\\
        \dots & \dots & \dots & \dots & \dots & \dots & \dots & \dots & \dots\\
        0 & \dots & \dots & \dots & a_k & 0 & \dots & \dots & -b_{k-1}\\
        0 & b_2 & b_4 & \dots & b_{2k-2} & a_1 & a_3 & \dots & a_{2k-3}\\
        0 & b_1 & b_3 & \dots & b_{2k-3} & 1 & a_2 & \dots & a_{2k-4}\\
        \dots & \dots & \dots & \dots & \dots & \dots & \dots & \dots & \dots\\
        0 & \dots & \dots & \dots & b_k & 0 & \dots & \dots & a_{k-1}\end{matrix}\right\vert,
    \end{align*}
    $k=2,\dots,n$ with $a_r = 0 = b_r$ whenever $r>n$.
\end{thm}
With the sake of keeping the paper pedagogical and focused on Algorithm 1 and its practical usage, a sketch of the proof of Theorem \ref{Thm_Cont_Frac} is given in Appendix \ref{app:cont_frac}. To make the connection between the determinants $\Delta_k$ and the coefficients computed via Algorithm 1, it can be shown that the following relations hold
\begin{align*}
    a_1^{(1)} &= \Delta_1,\\
    a_2^{(2)} &= \Delta_2,\\
    a_3^{(3)} &= \Delta_1^2\Delta_3,\\
    a_4^{(4)} &= \Delta_1^4\Delta_2^2\Delta_4,\\
    \dots,
\end{align*}
meaning that any $a_k^{(k)}, k=1,\dots,n$ may be expressed as $a_{k}^{(k)} = \Pi_{j=1}^k\Delta_j^{p_{jk}}$, where $p_{jk},j=1,\dots,k,$ are natural numbers that are not necessarily known a priori. Thanks to Theorem \ref{Thm_Determinants}, it is then easy to see that $a_{k}^{(k)} > 0$ if and only if $\Delta_k > 0, k=1,\dots,n$, which is a proof that Algorithm 1 is convergent.
}

\section{Example for a $4$-th degree polynomial}
\noindent
\label{sec:example}

Let us now show, as a pedagogical example, the table of coefficients for a $4$-th order polynomial of the form \eqref{eq:poly1} with $n=4$ and explicitly find the stability conditions.\\

The coefficients obtained with the method described in the previous section are given in Table \ref{tab:1}.
\begin{figure*}
\begin{center}
\begin{tabular}{ |p{3.5cm}|p{3.5cm}|p{3.5cm}|p{1.8cm}|}
\hline
${\color{black}a_1^{(1)}=~}a_1$  & ${\color{black}b_2^{(1)}=~}b_2$ & ${\color{black}a_3^{(1)}=~}a_3$ & ${\color{black}b_4^{(1)}=~}b_4$ \\
$b_1^{(1)}=a_1b_1-b_2$ & $a_2^{(1)}=a_1a_2-a_3$ & $b_3^{(1)}=a_1b_3-b_4$ & $a_4^{(1)}=a_1a_4$\\
\hline
$a_2^{(2)}=a_1a_2^{(1)}+b_1^{(1)}b_2$ & $b_3^{(2)}=a_1b_3^{(1)}-b_1^{(1)}a_3$ & $a_4^{(2)}=a_1a_4^{(1)}+b_1^{(1)}b_4$ &\\
$b_2^{(2)}=a_2^{(2)}b_2-a_1b_3^{(2)}$ & $a_3^{(2)}=a_2^{(2)}a_3-a_1a_4^{(2)}$ & $b_4^{(2)}=a_2^{(2)}b_4$ &\\
\hline
$a_3^{(3)}=a_2^{(2)}a_3^{(2)}+b_2^{(2)}b_3^{(2)}$ & $b_4^{(3)}=a_2^{(2)}b_4^{(2)}-b_2^{(2)}a_4^{(2)}$ & &\\
$b_3^{(3)}=a_3^{(3)}b_3^{(2)}-a_2^{(2)}b_4^{(3)}$ & $a_4^{(3)}=a_3^{(3)}a_4^{(2)}$ & &\\
\hline
$a_4^{(4)}=a_3^{(3)}a_4^{(3)}+b_3^{(3)}b_4^{(3)}$ & & &\\

 \hline
\end{tabular}\captionof{table}{Coefficients of the generalized Routh-Hurwitz criterion}\label{tab:1}
\end{center}
\end{figure*}

The necessary and sufficient conditions for the stability of the polynomial \eqref{eq:poly1} with $n=4$ are then given by 

\begin{equation}\label{eq:genRH}
\begin{cases}
a_{1}>0 \\
a_{2}^{(2)}=a_{1} a_{2}^{(1)}+b_1^{(1)}b_2>0 \\
a_{3}^{(3)}=a_{2}^{(2)} a_{3}^{(2)}+b_2^{(2)}b_3^{(2)}>0 \\
a_{4}^{(4)}=a_{3}^{(3)} a_{4}^{(3)}+b_3^{(3)}b_4^{(3)}>0
\end{cases}
\end{equation}

{\color{black}A straightforward computation}, which can be found in Appendix\ref{appendix:expl_derivation}, gives the expression for the generalized Routh-Hurwitz conditions 

\begin{equation}\label{eq:genRH_final}
\begin{cases}
a_{1}>0, \\
a_1^2a_2-a_1a_3-a_1b_1b_2+b_2^2:=\beta>0, \\
\beta^2a_3-\beta a_1^3a_4+\beta(a_1b_4+a_3b_2)(a_1b_1-b_2)\\
\hspace{0.5cm}-\beta a_1b_2(a_1b_3-b_4)+a_1[a_1(a_1b_3-b_4)\\
\hspace{0.5cm}-a_3(a_1b_1-b_2)]^2:=\gamma>0, \\
\gamma^2[a_1^2a_4-b_4(a_1b_1-b_2)]\\
\hspace{0.5cm}-\eta[\gamma a_1(a_1b_3-b_4)\gamma a_3(a_1b_1-b_2)-\beta\eta]>0,
\end{cases}
\end{equation} where \begin{align*}
    \eta&=\beta^2b_3-[\beta b_2-a_1(a_1(a_1b_3-b_4)-a_3(a_1b_1-b_2))]\varepsilon,
    \\
    \varepsilon &= [a_1^2a_4-b_4(a_1b_1-b_2)].
\end{align*}
From the above conditions \eqref{eq:genRH_final}, one can recover the classical Routh-Hurwitz conditions when the polynomial has real coefficients, i.e., $b_j=0,~\forall j\in\{1,2,3,4\}$, as we show in Appendix\ref{appendix:real_coeff}.

\section{Application to the PI-regulation of a rotating shaft}
\noindent
\label{sec:application}
An application arising from the theory of rotating shafts is considered in this section. We emphasize that this application has been chosen in order to highlight the efficiency of the generalized Routh-Hurwitz criterion, reason why some physical concepts are omitted. 
The dynamics of rotating shafts has triggered the interest of scholars since the very first studies carried out by William John Macquorn Rankine in 1869, with many notable results. For an overview on the topic, the reader may refer to \cite{Loewy_1970}. 
{\color{black}Therein, it is shown that such systems may be modeled by the following ordinary differential equation
\begin{equation}
    \ddot{x}(t) + (2k\omega + i 2\Omega )\dot{x}(t) + (\omega^2-\Omega^2)x(t) = u(t),
    \label{eq:RotatingShaft}
\end{equation}
where $k, \Omega$ and $\omega$ indicate a normalized dimensionless damping coefficient, an angular velocity $[s^{-1}]$ and the frequency of undamped oscillations $[s^{-1}]$, respectively. The quantity $u$ acts on the system as an external force. Before going to the control objective for that system, let us write it as a state-space model with state vector $\mathbf{x}(t) = (\begin{matrix} x(t) & \dot{x}(t)\end{matrix})^T$ and input function $u(t)$. There holds
\begin{equation}
    \mathbf{x}(t) = \mathbf{A}\mathbf{x}(t) + \mathbf{B}u(t), \mathbf{x}(0) = \mathbf{x}_0,
    \label{State-Space_Model}
\end{equation}
where $\mathbf{x}_0\in\mathbb{R}^2$ denotes the initial condition. In the state-space description \eqref{State-Space_Model}, the matrices $\mathbf{A}$ and $\mathbf{B}$ are given by
\begin{align*}
    \mathbf{A} = \left(\begin{matrix}0 & 1\\ \Omega^2-\omega^2 & -(2k\omega + i2\Omega)\end{matrix}\right), \mathbf{B} = \left(\begin{matrix}0 \\ 1\end{matrix}\right).
\end{align*}
As a control objective for \eqref{State-Space_Model}, let us consider the regulation of the position $x(t)$ to a constant prescribed reference position denoted by $x_r$. To achieve such an objective, one will rely on the well-established Proportional Integral (PI) control, see \cite{Borase} and references therein. In that way, the input $u(t)$ will take the following form
\begin{equation}
    u(t) = k_p(x(t)-x_r) + k_I \ell(t),
    \label{PI_Action}
\end{equation}
where the proportional and the integral gains denoted by $k_p$ and $k_I$ need to be determined for the closed-loop system to be stable. The quantity $\ell(t)$ is updated adaptively as
\begin{equation}
    \dot{\ell}(t) = x(t) - x_r.
    \label{Integral_Action}
\end{equation}
In that way, the closed-loop system composed of \eqref{State-Space_Model}, \eqref{PI_Action} and \eqref{Integral_Action} reads as
\begin{align}
         \left(\begin{smallmatrix}\dot{\mathbf{x}}(t)\\ \dot{\ell}(t)\end{smallmatrix}\right)
         &= \tilde{\mathbf{A}}\left(\begin{smallmatrix}\mathbf{x}(t)\\ \ell(t)\end{smallmatrix}\right) + \tilde{\mathbf{B}}x_r 
         \label{Closed-Loop_System}
\end{align}
where $\tilde{\mathbf{A}}=\left(\begin{smallmatrix}0 & 1 & 0\\ k_p+\Omega^2-\omega^2 & -(2k\omega + i 2\Omega ) & k_I\\ 1 & 0 & 0\end{smallmatrix}\right)$ and $\tilde{\mathbf{B}} = (\begin{matrix}0 & -k_p & -1\end{matrix})^T$.}
From \eqref{Closed-Loop_System}, if the gains $k_p$ and $k_I$ are chosen such that the matrix $\tilde{\mathbf{A}}$ is stable, then the control objective will be satisfied. In particular, the quantity $x$ will reach the equilibrium $x_r$. One needs then to determine in which cases the matrix $\tilde{\mathbf{A}}$ is a stable matrix. First observe that the characteristic polynomial of that matrix is given by 
\begin{align}
    q(s) &= \vert sI - \tilde{\mathbf{A}}\vert\nonumber\\
    &= s^3 + (2k\omega +i 2\Omega )s^2 + (\omega^2 - \Omega^2 - k_p)s - k_I.
    \label{Char_Poly_Appli}
\end{align}
We shall therefore rely on the generalized Routh-Hurwitz criterion detailed in Algorithm \ref{Alg:RH}. The consecutive arrays of numbers generated by this algorithm are given in Table \ref{tab:appli}.
\begin{figure*}
\begin{tabular}{|p{5.2cm}|p{5.2cm}|p{5.2cm}|}
\hline
$a_1^{(1)} = 2k\omega$ & $b_2^{(1)} = 0$ & $a_3^{(1)} = -k_I$\\
$b_1^{(1)} = 4k\omega\Omega$ & $a_2^{(1)} = 2k\omega(\omega^2-\Omega^2-k_p)+k_I$ & $b_3^{(1)} = 0$\\
\hline
$a_2^{(2)} = 2k\omega(2k\omega(\omega^2-\Omega^2-k_p)+k_I)$ & $b_3^{(2)} = k_I(4k\omega\Omega)$ & \\
$b_2^{(2)} = -8k_I k^2 \omega^2\Omega$ & $a_3^{(2)} = -k_I a_2^{(2)}$ & \\
\hline 
$a_3^{(3)} = -k_I(a_2^{(2)})^2 - 32 k_I^2 k^3\omega^3\Omega^2$ & & \\
\hline
\end{tabular}
\captionof{table}{Generalized Routh-Hurwitz table for the polynomial \eqref{Char_Poly_Appli}}\label{tab:appli}
\end{figure*}
According to Algorithm \ref{Alg:RH}, the matrix $\tilde{\mathbf{A}}$ is stable if and only if the following three conditions are satisfied
\begin{align}
2k\omega &> 0\label{cond_1_Appli}\\
2k\omega\left[2k\omega(\omega^2 - \Omega^2-k_p) + k_I\right]&>0\label{cond_2_Appli}\\
-8k_I^2k\omega\Omega^2 - k_I\left[2k\omega(\omega^2 - \Omega^2-k_p) + k_I\right]^2&>0\label{cond_3_Appli}.
\end{align}
In order to illustrate the feasibility of Conditions \eqref{cond_1_Appli}--\eqref{cond_3_Appli}, a grid has been made with different values of $k_I$ and $k_p$. At each point of the grid, Conditions \eqref{cond_1_Appli}--\eqref{cond_3_Appli} have been tested. If these conditions are all satisfied, the value $1$ (yellow part) has been placed on the grid and $0$ (blue part) otherwise. The resulting picture is depicted\footnote{{\color{black}The software used for all the Figures in this paper is \cite{MATLAB2021}, version R2021a.}} in Figure \ref{fig1:Colormap1}. As a matter of comparison, another test has been performed. For each point of the grid $(k_I,k_p)$, the value of the largest real part of the eigenvalues of the matrix $\tilde{\mathbf{A}}$ has been encoded in the grid. A contour plot has then been performed and it is shown in Figure \ref{fig2:contourPlot}. From the latter, it is clear that the stability region is the same as the one highlighted with Figure \ref{fig1:Colormap1}, obtained from the stability conditions \eqref{cond_1_Appli}--\eqref{cond_3_Appli}. In the above derivation, the parameters have been set to $k=1, \omega=2$ and $\Omega=2$.

\begin{figure}
    \centering
    \includegraphics[scale=0.55,trim=3cm 2cm 1cm 2cm]{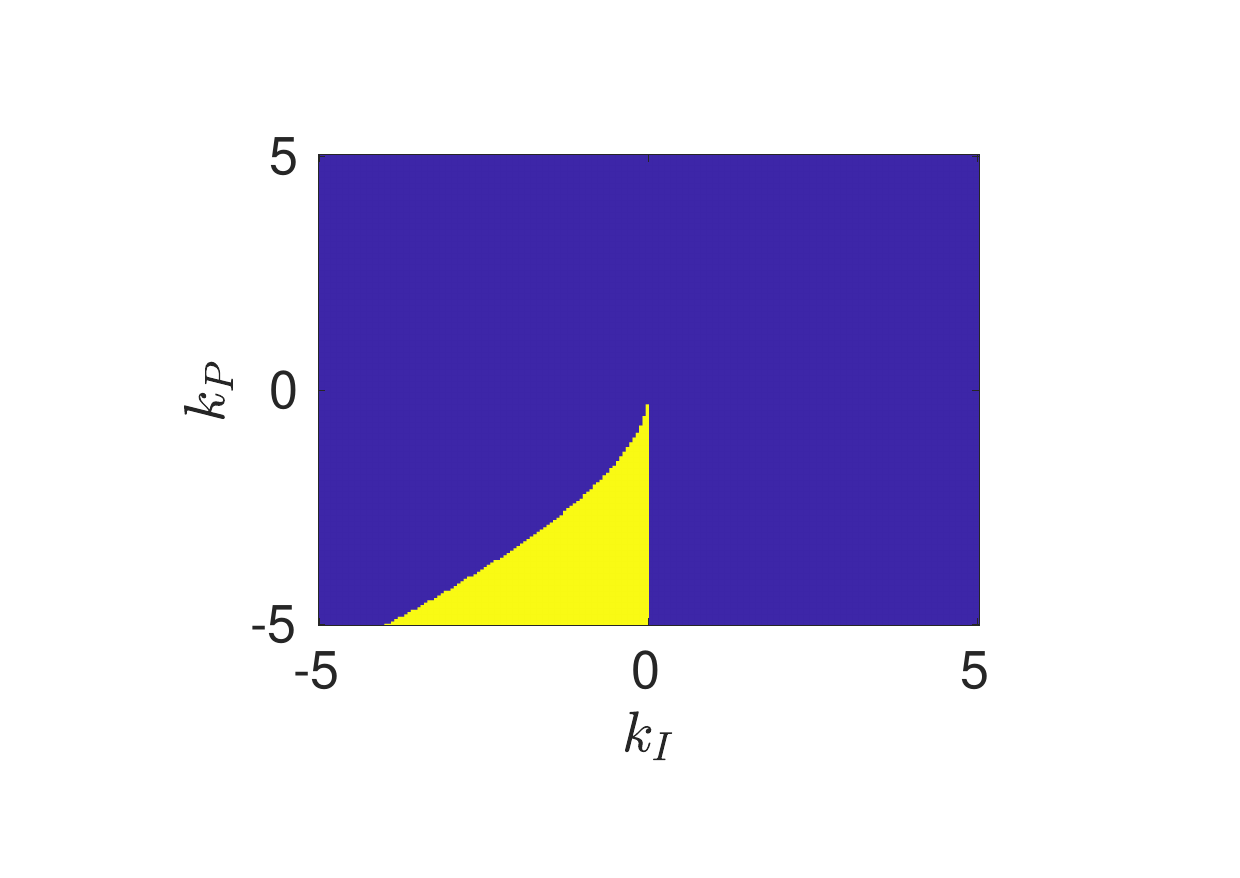}
    \caption{A grid $(k_I,k_p)$ at which each point is either $1$ or $0$, based on conditions \eqref{cond_1_Appli}--\eqref{cond_3_Appli}.}
    \label{fig1:Colormap1}
\end{figure}

\begin{figure}
    \centering
    \includegraphics[scale=0.55,trim=3cm 2cm 1cm 2cm]{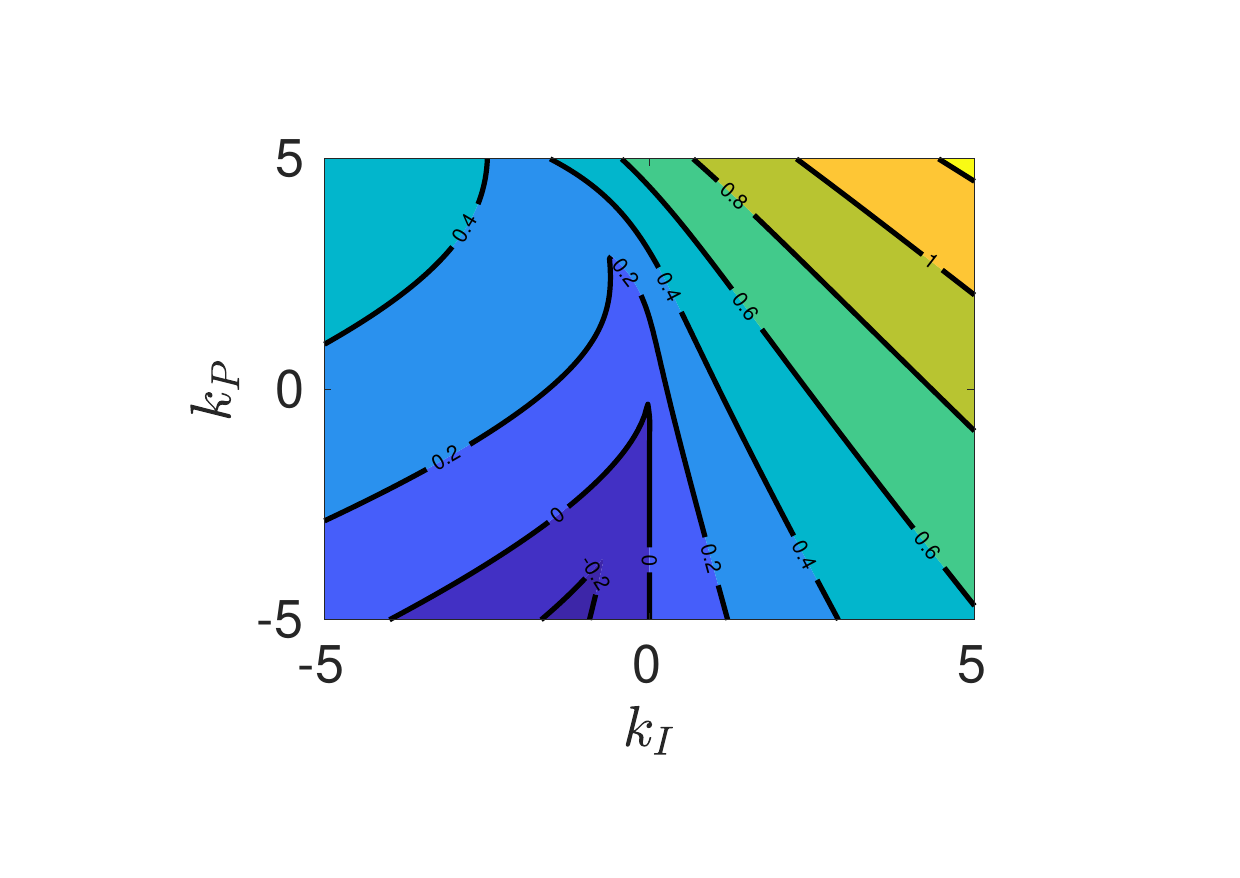}
    \caption{Contour plot of the function $\lambda_s(k_I,k_p) := \max\{\mathfrak{Re}(\lambda), \lambda\in\sigma(\tilde{\mathbf{A}})\}$}
    \label{fig2:contourPlot}
\end{figure}

{\color{black}In order to illustrate the efficiency of the PI control action, the system responses are depicted in Figures \ref{fig:x(t)}, \ref{fig:dx(t)} and \ref{fig:ell(t)} for the following values of $x_r, k_I$ and $k_p$: $x_r = 2, k_I = -1.18, k_p = -3.59$. Note that $k_I$ and $k_p$ have been chosen in accordance to the stability region of the matrix $\tilde{\mathbf{A}}$. In these responses, one can observe that the trajectory $x(t)$ converges to the reference signal $x_r$ while both $\dot{x}(t)$ and $\ell(t)$ tend to $0$ as $t$ goes to $\infty$.
\begin{figure}
    \centering
    \includegraphics[scale=0.55,trim=3cm 2cm 1cm 2cm]{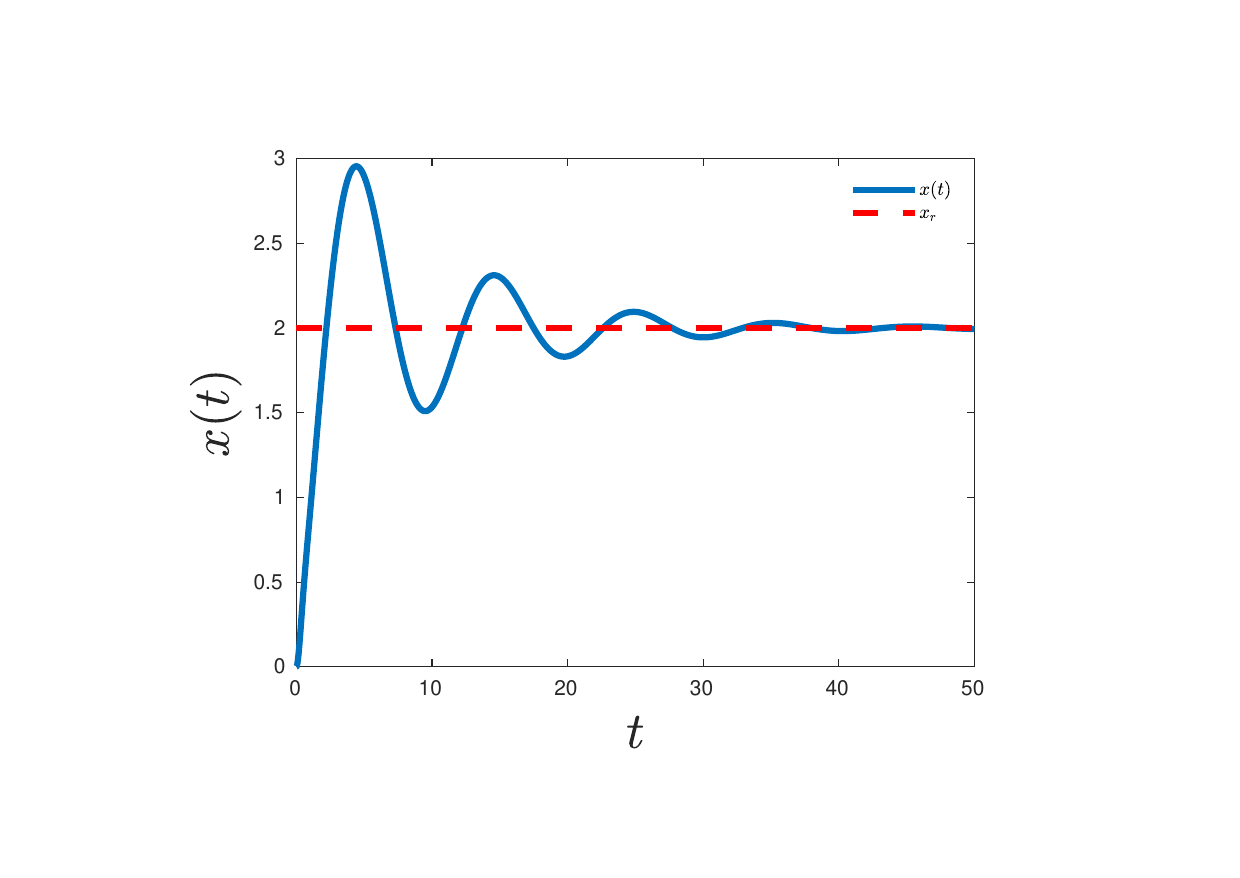}
    \caption{Trajectory $x(t)$ for $x_r = 2, k_I = -1.18, k_p = -3.59$.}
    \label{fig:x(t)}
\end{figure}
\begin{figure}
    \centering
    \includegraphics[scale=0.55,trim=3cm 2cm 1cm 2cm]{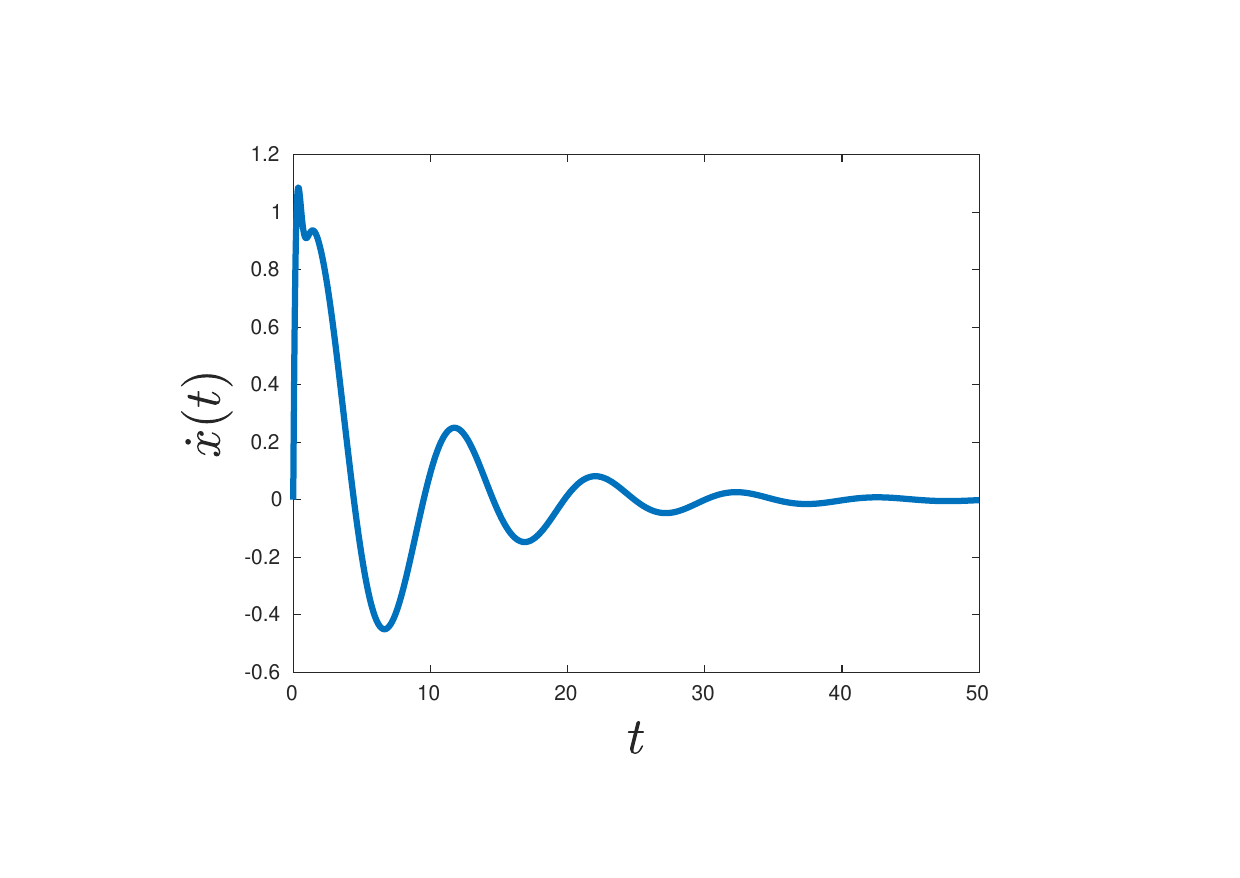}
    \caption{Trajectory $\dot{x}(t)$ for $x_r = 2, k_I = -1.18, k_p = -3.59$.}
    \label{fig:dx(t)}
\end{figure}
\begin{figure}
    \centering
    \includegraphics[scale=0.55,trim=3cm 2cm 1cm 2cm]{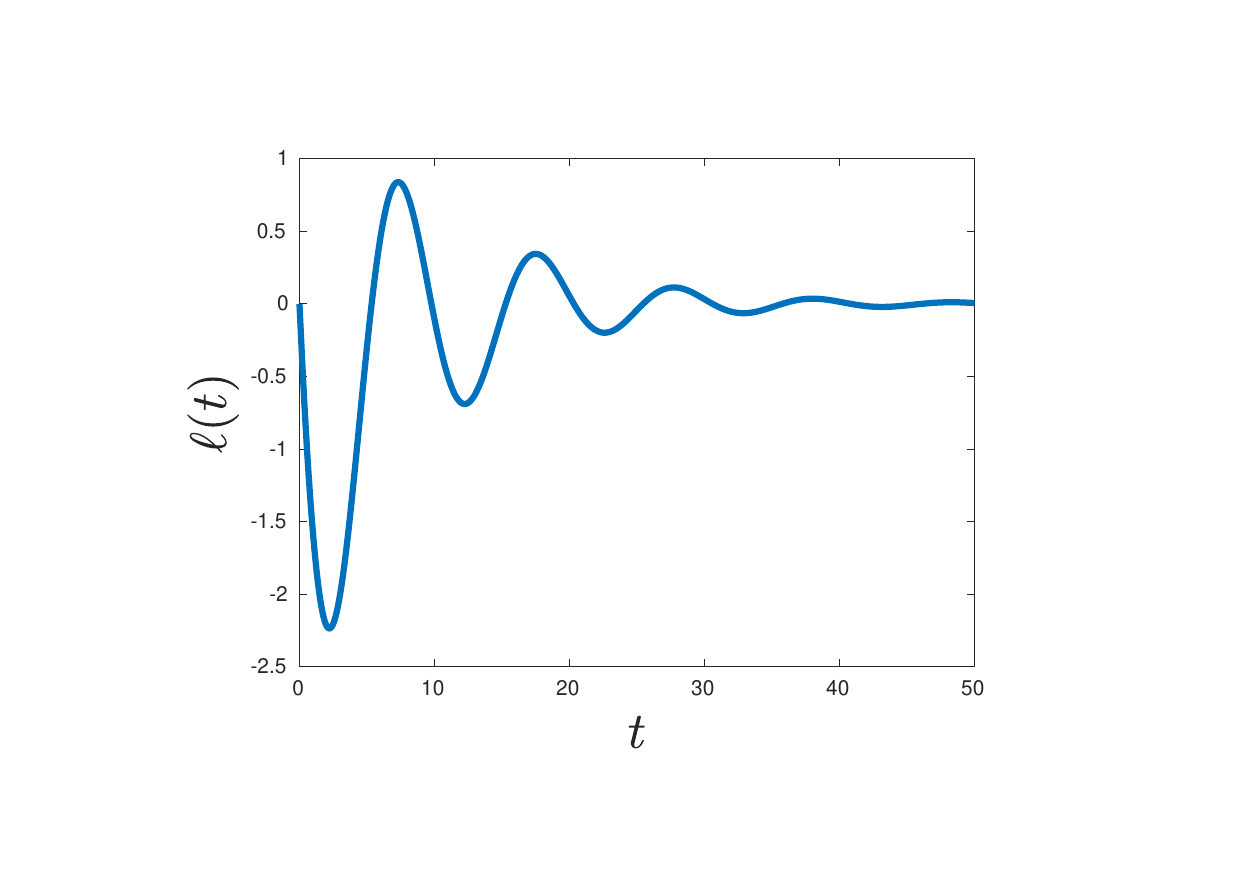}
    \caption{Trajectory $\ell(t)$ for $x_r = 2, k_I = -1.18, k_p = -3.59$.}
    \label{fig:ell(t)}
\end{figure}
}

{\color{black}
\begin{rem}
An interesting perspective regarding the presented example would be the study of the long-term behavior of the system after being exposed to aging and to some noise. In particular, questions like "How the gains in the PI-controller would be affected in such a situation?", "How would be the consequences on the stability of the system?" should be investigated. References like \cite{Wang_2021} and \cite{Dai_2021} could be a starting point. Therein, reliability of dynamical systems is studied thanks to different notions such as for instance the \textit{principle of maximum entropy}, ... . These results are applied to systems composed of electrical circuits or dynamical systems subject to wear and vibration.
\end{rem}
}


\section{Conclusion}
\noindent
\label{sec:ccl}
In this note, we have clarified and explained in a constructive and pedagogical way an extension of the classical Routh-Hurwitz criterion to polynomials with complex coefficients. The general algorithm to determine whether the roots of such a polynomial are located in $\mathbb{C}_0^-$ or not is given in Section \ref{sec:method}, {\color{black} broken down in a pedagogical way}. Then, {\color{black} the latter is explicitly derived} for a $4$-th order polynomial in Section \ref{sec:example}. Finally, an application to the PI regulation of a rotating shaft whose own dynamics exhibit complex coefficients is detailed in Section \ref{sec:application}, giving rise to the study of the stability of a $3$-rd order polynomial. Our presentation of the algorithm and the given examples make the method understandable and ready to use also for scholars and students outside the control community. {\color{black} Our work paves} the way for further advancements in applications where complex polynomials appear{\color{black}, such as dynamics on networks and hypergraphs, where asymmetric topologies lead to complex coefficients.}

\section*{Acknowledgements} The authors are grateful to Alice Bellière for useful discussions and feedback {\color{black} and to two anonymous Reveiwers, whose comments and constructive criticism have improved the quality of this work}. This research was conducted with the financial support of F.R.S-FNRS. A.H. is supported by a FNRS Postdoctoral Fellowship, Grant CR 40010909. During the realization of this work, R.M. was supported by a FRIA Fellowship, funded by the Walloon Region, Grant FC 33443.

\appendix

\section{Explicit derivation of the Generalized Routh-Hurwitz conditions}\label{appendix:expl_derivation}

Let us now explicitly derive the conditions \eqref{eq:genRH_final}.\\

Derivation of the second condition:
$$
\begin{aligned}
a_{2}^{(2)} &=a_{1} a_{2}^{(2)}+b_{1}^{(1)} b_{2}  =a_{1}(a_{1} a_{2}-a_{3})+b_{2}(a_{1} b_{1}-b_{2})\\
& =a_{1}^{2} a_{2}-a_{1} a_{3}+a_{1} b_{1} b_{2}-b_{2}^{2} := \beta \\
\end{aligned}
$$

Derivation of the third condition:

\begin{align*}
a_{3}^{(3)}&=a_{2}^{(2)} a_{3}^{(2)}+b_{2}^{(2)} b_{3}^{(2)}\\
&=a_{2}^{(2)}(a_{2}^{(2)} a_{3}-a_{1} a_{4}^{(2)})\\
&+(a_{2}^{(2)} b_{2}-a_{1} b_{3}^{(2)})(a_{1}b_3^{(1)}-b_{1}^{(1)} a_{3}) \\
& =(a_{2}^{(2)})^2a_{3}-a_{1} a_{2}^{(2)} a_{4}^{(2)}+a_{1} b_{2} a_{2}^{(2)} b_{3}^{(1)}-a_{3} b_{2} a_{2}^{(2)} b_{1}^{(1)}\\
&-a_{1}^{2} b_{3}^{(1)} b_{3}^{(2)}+a_{1} a_{3} b_{1}^{(1)} b_{3}^{(2)}
\end{align*}

\begin{align*}
& =\beta^{2} a_{3}-\beta a_{1} (a_{1} a_{4}^{(1)}+b_{1}^{(1)} b_{4})+\beta a_{1} b_{2} (a_{1} b_{3}-b_{4})\\
&-\beta a_{3} b_{2} (a_{1} b_{1}-b_{2})-a_{1}^{2}\left(a_{1} b_{3}-b_{4}\right)\left(a_{1} b_{3}^{(1)}-a_{3} b_{1}^{(1)}\right)\\
&+a_{1} a_{3}\left(a_{1} b_{1}-b_{2}\right)\left(a_{1} b_{3}^{(1)}-b_{1}^{(1)} a_{3}\right)\\
&=\beta^{2} a_{3}-\beta a_{1} \left[a_{1}^{2} a_{4}+b_{4}\left(a_{1} b_{1}-b_{2}\right)\right]\\
&+\beta a_{1} b_{2} \left(a_{1} b_{3}-b_{4}\right)-\beta a_{3} b_{2} \left(a_{1} b_{1}-b_{2}\right)\\
& -a_{1}^{2}\left(a_{1} b_{3}-b_{4}\right)\left[a_{1}\left(a_{1} b_{3}-b_{4}\right)-a_{3}\left(a_{1} b_{1}-b_{2}\right)\right]\\
&+a_{1} a_{3}\left(a_{1} b_{1}-b_{2}\right)\left[a_{1}\left(a_{1} b_{3}-b_{4}\right)-a_{3}\left(a_{1} b_{1}-b_{2}\right)\right]\\
& =\beta^{2} a_{3}-\beta a_{1} \left(a_{1}^{2} a_{4}+b_{4}\left(a_{1} b_{1}-b_{2}\right)\right)\\
&+\beta  a_{1} b_{2} \left(a_{1} b_{3}-b_{4}\right) - \beta a_{3} b_{2} \left(a_{1} b_{1}-b_{2}\right)\\
& -\left[a_{1}^{2}\left(a_{1} b_{3}-b_{4}\right)-a_{1} a_{3}\left(a_{1} b_{1}-b_{2}\right)\right]\dots\\
&\dots\left[a_{1}\left(a_{1} b_{3}-b_{4}\right)-a_{3}\left(a_{1} b_{1}-b_{2}\right)\right]\\
\end{align*}

$$\begin{aligned}
&=\beta^{2} a_{3}-\beta a_{1}^{3} a_{4}-\beta  a_{1} b_4\left(a_{1} b_{1}-b_{2}\right)+\beta a_{1} b_{2} \left(a_{1} b_{3}-b_{4}\right)\\
&-\beta a_{3} b_{2} \left(a_{1} b_{1}-b_{2}\right)-a_{1}\left[a_{1}\left(a_{1} b_{3}-b_{4}\right)-a_{3}\left(a_{1} b_{1}-b_{2}\right)\right]^{2}\\
&=\beta^2a_3-\beta a_1^3a_4-\beta(a_1b_4+a_3b_2)(a_1b_1-b_2)\\
&+\beta a_1b_2(a_1b_3-b_4)-a_1[a_1(a_1b_3-b_4)-a_3(a_1b_1-b_2)]^2\\
&:=\gamma
\end{aligned}
$$

Derivation of the fourth condition:

$$
\begin{aligned}
a_{4}^{(4)}&=\gamma a_{4}^{(3)}+b_{3}^{(3)} b_{4}^{(3)} \\
& =\gamma\left(\gamma a_{4}^{(2)}\right)+\left(\gamma b_{3}^{(2)}-\beta b_{4}^{(3)}\right) b_{4}^{(3)} \\
&=\gamma^{2}\left(a_{1} a_{4}^{(1)}+b_{1}^{(1)} b_{4}\right)\\
&+\gamma\left(a_{1} b_{3}^{(1)}-a_{3} b_{1}^{(1)}\right)\left(\beta b_{4}^{(2)}-b_{2}^{(2)} a_{4}^{(2)}\right)-\beta\left(b_{4}^{(3)}\right)^{2}\\
& =\gamma^{2}\left[a_{1}^{2} a_{4}+b_{4}\left(a_{1} b_{1}-b_{2}\right)\right]\\
&+\gamma\left[a_{1}\left(a_{1} b_{3}-b_{4}\right)-a_{3}\left(a_{1} b_{1}-b_{2}\right)\right] \eta - \beta (b_4^{(3)})^2 
\end{aligned}
$$
Where we have defined
$$
\begin{aligned}
 \eta&:=\beta^{2} b_{4}-\left(\beta b_{2}-a_{1} b_{3}^{(2)}\right)\left(a_{1} a_{4}^{(1)}+b_{1}^{(1)} b_{4}\right)\\
 & = \beta^{2} b_{4}-\left[\beta b_{2}-a_{1}\left(a_{1} b_{3}^{(1)}-a_{3} b_{1}^{(1)}\right)\right]\varepsilon,
\end{aligned}
$$
with $\varepsilon := \left[a_{1}^{2} a_{4}+b_{4}\left(a_{1} b_{1}-b_{2}\right)\right]$. Moreover, we can rewrite the coefficient $b_4^{(3)}$ in terms of $\eta$, obtaining
\begin{align*}
 b_4^{(3)}&= \beta^{2} b_{4}-b_{2}^{(2)} a_{4}^{(2)}\\
 &= \beta^{2} b_4-\left[\beta b_{2}-a_{1}\left(a_{1} b_{3}^{(1)}-a_{3} b_{1}^{(1)}\right)\right]\varepsilon\\
 &=\eta
\end{align*}

Hence, the explicit expression of the fourth condition becomes

\begin{align*}
    a_{4}^{(4)}&=\gamma^2\varepsilon+\eta[\gamma a_1(a_1b_3-b_4)-\gamma a_3(a_1b_1-b_2)-\beta\eta]
\end{align*}

Lastly, let us derive the explicit expression for $\eta$:
$$
\begin{aligned}
 \eta&:=\beta^2b_4-[\beta b_2-a_1(a_1(a_1b_3-b_4)-a_3(a_1b_1-b_2))]\varepsilon.
\end{aligned}
$$
\section{Attainment of Classical Routh-Hurwitz criterion in case of real coefficient}\label{appendix:real_coeff}

Let us now show that the stability conditions are the same as the classical Routh-Hurwitz criterion in case of real coefficients, namely $b_{j}=0, \forall j\in\{0,..,4\}$. For simplicity, let us again consider a $4$-order polynomial

\begin{displaymath}
    p(s)=s^4+a_1s^3+a_2s^2+a_3 s+a_4
\end{displaymath}

The table of the coefficients is given by
\begin{center}
\begin{tabular}{ |p{0.5cm}||p{3.3cm}|p{0.5cm}|p{0.5cm}|p{0.5cm}| }
 
 \hline
     & & & & \\
$s^4$& $1$ & $a_2$ & $a_4$ & $0$ \\
 \hline
   & & & &  \\
 $s^3$   & $a_1$ & $a_3$ & $0$ & $0$ \\
 \hline
   & & & & \\
$s^2$   & $\displaystyle\frac{a_1a_2-a_3}{a_1}$  & $a_4$ & $0$ & $0$ \\
\hline
  & & & & \\
 $s^1$ & $\displaystyle\frac{(a_1a_2-a_3)a_3-a_1^2a_4}{a_1a_2-a_3}$ & $0$ & $0$ & $0$ \\
 \hline
   & & & & \\
 $s^0$    & $a_4$ &$0$ & $0$&$0$  \\

 \hline
\end{tabular}\captionof{table}{Coefficients of the classical Routh-Hurwitz criterion}\label{tab:2}
\end{center}

and the necessary and sufficient stability conditions are given by

\begin{equation}\label{eq:classicRH}
    \begin{cases}
 a_{1}>0 \\ a_{1} a_{2}-a_{3}>0 \\ a_{4}>0 \\ \left(a_{1} a_{2}-a_{3}\right) a_{3}-a_{1}^{2} a_{4}>0 
    \end{cases}
\end{equation}

When $b_j=0~\forall j\in\{0,..,4\}$, Table \ref{tab:1} becomes the table that is given in \ref{tab:3}.

\begin{figure*}
\begin{center}
\begin{tabular}{ |p{2.4cm}|p{3.2cm}|p{2cm}|p{1.8cm}|}
\hline
$a_1$  & $0$ & $a_3$ & $0$\\
$b_1^{(1)}=0$ & $a_2^{(1)}=a_1a_2-a_3$ & $b_3^{(1)}=0$ & $a_4^{(1)}=a_1a_4$\\
\hline
$a_2^{(2)}=a_1a_2^{(1)}$ & $b_3^{(2)}=0$ & $a_4^{(2)}=a_1a_4^{(1)}$ &\\
$b_2^{(2)}=0$ & $a_3^{(2)}=a_2^{(2)}a_3-a_1a_4^{(2)}$ & $b_4^{(2)}=0$ &\\
\hline
$a_3^{(3)}=a_2^{(2)}a_3^{(2)}$ & $b_4^{(3)}=0$ & &\\
$b_3^{(3)}=0$ & $a_4^{(3)}=a_3^{(3)}a_4^{(2)}$ & &\\
\hline
$a_4^{(4)}=a_3^{(3)}a_4^{(3)}$ & & &\\
\hline
\end{tabular}\captionof{table}{Coefficients of the generalized Routh-Hurwitz criterion for $b_j=0$}\label{tab:3}
\end{center}
\end{figure*}
and the necessary and sufficient conditions for stability are given by Eq. \eqref{eq:genRH}, which for real coefficients become

\begin{equation}\label{eq:class_recovered}
\begin{cases}
a_{1}>0 \\
a_{2}^{(2)}=a_{1} a_{2}^{(1)}>0 \\
a_{3}^{(3)}=a_{2}^{(2)} a_{3}^{(2)}>0 \\
a_{4}^{(4)}=a_{3}^{(3)} a_{4}^{(3)}>0
\end{cases}
\end{equation}

Let us remember that they all have to stand simultaneously, in order for the system to be stable. The equivalence between the first condition of Eq. \eqref{eq:class_recovered} and the first of \eqref{eq:classicRH} is trivial. The second condition of Eq. \eqref{eq:class_recovered} gives us \begin{displaymath}
    a_{2}^{(2)}=a_{1} a_{2}^{(1)}=a_1( a_{1} a_{2}-a_{3})>0
\end{displaymath} Since $a_1>0$, we have the second condition of Eq. \eqref{eq:classicRH}. The third condition of Eq. \eqref{eq:class_recovered} gives us 
\begin{displaymath}
   a_3^{(3)} =a_{1}\left(a_{1} a_{2}-a_{3}\right)\left[a_{1}\left(a_{1} a_{2}-a_{3}\right) a_{3}-a_{1}^{3} a_{4}\right]>0
\end{displaymath} Given the first two conditions, namely $a_1>0$ and $a_{1} a_{2}-a_{3}>0$, we obtain $\left(a_{1} a_{2}-a_{3}\right) a_{3}-a_{1}^{3} a_{4}>0$, which is exactly the fourth condition of Eq. \eqref{eq:classicRH}. Lastly, the fourth condition of Eq. \eqref{eq:class_recovered} gives us \begin{displaymath}
    a_{4}^{(4)}=(a_3^{(3)})^2a_1^2a_4>0
\end{displaymath} 
which reduces to $a_4>0$, i.e., the third condition of Eq. \eqref{eq:classicRH}. \\
Hence, from the generalized Routh-Hurwitz conditions for the case of real coefficients, we attained the classical Routh-Hurwitz conditions, proving the equivalence.

\section{Comparison with an existing method}\label{app:comparison}

{\color{black} In this Appendix we compare the method hereby developed with the method developed in \cite{carletti} for a $3$-th order polynomial \begin{equation}\label{eq:app_t3rd}
    q(s)=s^3+(a_1+ib_1)s^2+(a_2+ib_2)s+(a_3+ib_3).
    \end{equation}

The generalized Routh-Hurwitz criterion developed in the Main Text gives the following necessary and sufficient condition for the stability of the above polynomial: 
\begin{displaymath}
\begin{cases}
    a_1^{(1)}>0, \\
    a_2^{(2)}>0, \\
    a_3^{(3)}>0,
    \end{cases}
\end{displaymath} whose expression can be computed through the Algorithm \ref{Alg:RH}. Explicitly, we have:

 \begin{equation}\label{eq:rh_new_app_3rd}
\begin{cases}
 a_1>0, \\
  \nu := a_1^2a_2 - a_1 a_3 + a_1 b_1 b_2 - b_2^2>0, \\
 a_3\nu^2 + \nu [a_1^2 b_2 b_3 - a_3 b_2(a_1 b_1-b_2)] - a_1^5 b_3^2\\
 +2a_1^3 a_3 b_3(a_1 b_1-b_2) - a_1 a_3^2(a_1 b_1-b_2)^2>0.
   \end{cases}
   \end{equation}
\vspace{0.4cm}

The method of \cite{carletti} consists in multiplying $q(s)$ by the polynomial $\bar{q}(s)$, whose coefficients are the complex conjugate of the former. The authors prove that the obtained polynomial $Q(s)=q(s)\bar{q}(s)$, whose degree is double the degree of $q(s)$, has real coefficients with the same real parts of the coefficients of $q(s)$. Hence, the stability of $q(s)$ can be determined by applying the classic Routh-Hurwitz criterion to the polynomial $Q(s)$. \\

For our case of $3$-rd degree polynomial, we need to multiply Eq. \eqref{eq:rh_new_app_3rd} by \begin{displaymath}
    \bar{q}(s)=s^3+(a_1-ib_1)s^2+(a_2-ib_2)s^1+(a_3-ib_3), 
    \end{displaymath} obtaining thus the $6$-degree polynomial: \begin{equation}\label{eq:app_bigQ}
    Q(s)=s^{6}+A_1s^5+A_2s^4+A_3s^3+A_4s^2+A_5s+A_6,
 \end{equation} whose coefficients are given by: \begin{eqnarray*} 
    A_1&=& 2a_1,  \nonumber \\
    A_2&=& 2a_2+a_1^2+b_1^2, \nonumber \\
    A_3&=& 2a_3+2(a_1a_2+b_1b_2), \nonumber \\
    A_4&=& a_2^2+b_2^2+2(a_1a_3+b_1b_3), \nonumber \\
    A_5&=&  2(a_2a_3+b_2b_3),\nonumber \\
    A_6&=& a_3^2+b_3^2. \nonumber \\
\end{eqnarray*}

By applying the classical Routh-Hurwitz criterion, we obtain that polynomial $Q(s)$, and thus $q(s)$, is stable if and only if: \begin{equation}\label{eq:app_ugly}
\begin{cases}
\scriptstyle A_1>0, \\
 \scriptstyle A_1A_2-A_3>0, \\
 \scriptstyle  (A_1A_2-A_3)A_3-A_1(A_1A_4-A_5)>0, \\
 \scriptstyle \left\vert\begin{smallmatrix}A_1 & 1 & 0 & 0\\ A_3 & A_2 & A_1 & 1\\ A_5 & A_4 & A_3 & A_2\\ 0 & A_6 & A_5 & A_4\end{smallmatrix}\right\vert > 0,\vspace{0.2cm}\\
 \scriptstyle \left\vert\begin{smallmatrix}A_1 & 1 & 0 & 0 & 0\\ A_3 & A_2 & A_1 & 1 & 0\\ A_5 & A_4 & A_3 & A_2 & A_1\\ 0 & A_6 & A_5 & A_4 & A_3\\ 0 & 0 & 0 & A_6 & A_5\end{smallmatrix}\right\vert > 0,\vspace{0.2cm}\\
 \scriptstyle \left\vert\begin{smallmatrix}A_1 & 1 & 0 & 0 & 0 & 0\\ A_3 & A_2 & A_1 & 1 & 0 & 0\\ A_5 & A_4 & A_3 & A_2 & A_1 & 1\\ 0 & A_6 & A_5 & A_4 & A_3 & A_2\\ 0 & 0 & 0 & A_6 & A_5 & A_4\\ 0 & 0 & 0 & 0 & 0 & A_6\end{smallmatrix}\right\vert > 0.
   \end{cases}
   \end{equation}


We do not need to substitute the explicit expressions of the coefficients $A_i$ to see how much more complicated it is to compute the above conditions \eqref{eq:app_ugly} rather than the conditions \eqref{eq:rh_new_app_3rd} obtained with our method.

}

{\color{black}
\section{Sketch of the proof of Theorem \ref{Thm_Cont_Frac}}\label{app:cont_frac}
Let us first suppose that the continued fraction expansion \eqref{Continued_Frac} holds with $c_i > 0, i=0,\dots,n-1$. We then may rewrite it as generated by the following sequence
\begin{align}
    f &= \frac{c_0}{s+c_0+d_1+f_1}, f_1 = \frac{c_1}{s+d_2+f_2}, \dots,\nonumber\\
    f_{n-1} &= \frac{c_{n-1}}{s+d_n+f_n}, f_n = 0.
    \label{sequence_Cont_Frac}
\end{align}
As it is shown in \cite{Wall}, the relation $\vert f-\frac{1}{2}\vert\leq \frac{1}{2}$ holds when $\mathfrak{Re}(f_1)\geq 0$ and $\mathfrak{Re}(s)\geq 0$. Moreover, since $\mathfrak{Re}(f_p)\geq 0$ implies that $\mathfrak{Re}(f_{p-1})\geq 0$ for $p=2,\dots,n$, then 
\begin{equation*}
    \left\vert\frac{\mathfrak{p}(s)}{p(s)}-\frac{1}{2}\right\vert\leq\frac{1}{2}, \text{ for }\mathfrak{Re}(s)\geq 0.
\end{equation*}
As a direct consequence, the polynomial $p$ may not have any roots when $\mathfrak{Re}(s)\geq 0$, meaning that it is stable.
In a second time, let us assume that the polynomial $p$ is stable. By $\overline{p}(s)$, we denote the polynomial $p$ whose coefficients are complex conjugated. Then, one may observe that the polynomial $\mathfrak{p}(s)$ in \eqref{Aux_Poly} is either $\frac{1}{2}(p(s) + \overline{p}(-s))$ if $n$ is odd or $\frac{1}{2}(p(s) - \overline{p}(-s))$ if $n$ is even. As the roots of $p(s)$ and $\overline{p}(-s)$ are symmetrical to the imaginary axis, the roots of $\mathfrak{p}(s)$ lie on the imaginary axis, see \cite{Wall} for more detailed arguments. It then follows that the fraction $\mathfrak{p}(s)/p(s)$ is irreducible. Hence, it may be written as 
\begin{equation}
    \frac{\mathfrak{p}(s)}{p(s)} = \frac{c_0}{s + c_0 + d_1 + \frac{c(s)}{\mathfrak{p}(s)}},
    \label{Intermediate_Fraction}
\end{equation}
where $c_0$ is the opposite of the sum of the real parts of the roots of $p(s)$ (which is a positive number), $d_1$ is pure imaginary or zero, and $c(s)/\mathfrak{p}(s)$ is an irreducible rational fraction in which $c(s)$ is of degree less than $n-1$. Following \cite[Proof of Theorem A]{Wall}, one has that $\mathfrak{Re}(c(s)/\mathfrak{p}(s))\geq 0$ for $\mathfrak{Re}(s)\geq 0$. This has the consequence that 
\begin{equation*}
    \frac{c(s)}{\mathfrak{p}(s)} = \sum_{j=1}^{n-1} \frac{l_j}{s + ix_j},
\end{equation*}
where $x_j$ are real and distinct while $l_j>0, j=1,\dots,n-1$. It then follows that $\frac{-ic(-is)}{\mathfrak{p}(-is)} = \sum_{j=1}^{n-1} \frac{l_j}{s - x_j}$, which implies that
\begin{align*}
    \frac{-ic(-is)}{\mathfrak{p}(-is)} = \frac{c_1}{s + id_2 - \displaystyle\frac{c_2}{s + i d_3 - \displaystyle\frac{c_3}{s + i d_4 + \stackunder{}{\ddots\stackunder{}{\displaystyle{}- \frac{c_{n-1}}{\displaystyle
        s+ i d_n}}}}}},
\end{align*}
where $c_j, j=1,\dots,n-1$ are real and positive and $d_j, j=2,\dots,n$ are pure imaginary or zero. Replacing $s$ by $is$, dividing both sides by $-i$ and incorporating the obtained continued fraction into \eqref{Intermediate_Fraction} concludes the proof.
}

\section*{References}
\bibliographystyle{elsarticle-harv}     
\bibliography{biblio}


\end{document}